\documentclass[11pt]{amsart}
\usepackage{amsmath}
\usepackage{amscd}
\usepackage{amsfonts}
\usepackage{amssymb}
\usepackage{latexsym}
\usepackage{amscd}

\numberwithin{equation}{section}

\newtheorem{defn}{Definition}[section]
\newtheorem{theorem}{Theorem}[section]
\newtheorem{assumption}[theorem]{Assumption}

\newtheorem{corollary}[theorem]{Corollary}

\newtheorem{lemma}[theorem]{Lemma}
\newtheorem{prop}[theorem]{Proposition}
\newtheorem{remark}[theorem]{Remark}

\def \n{\noindent}

\def \v{\vskip 0.1in}

\def \mc{\mathcal}
\def \mf{\mathbf}

\def \cplane{\mathbb{C}}

\def \pone{\mathbb{P}}
\def \integer{\mathbb{Z}}

\def \mfp{\mathfrak{p}^s}
\def \mfq{\mathfrak{q}^s}
\def \mfP{\mathfrak{p}^{sf}}
\def \mfQ{\mathfrak{q}^{sf}}

\def \om{\overline{\mathcal{M}}}
\def \M{\mathcal{M}}

\def \mbx{\mathbf{x}}

\begin{document}

\title{Ruan's Conjecture on Singular symplectic flops}
\author{Bohui Chen}
\address{Department of Mathematics, Sichuan University,
        Chengdu,610064, China}
\email{bohui@cs.wisc.edu}
\author{An-Min Li}
\address{Department of Mathematics, Sichuan University,
        Chengdu,610064, China}
\email{math$\_$li@yahoo.com.cn}
\author{Guosong Zhao }
\address{Department of Mathematics, Sichuan University,
        Chengdu,610064, China}
\email{gszhao@scu.edu.cn}
\thanks{B.C. and A.L. are supported by  NSFC,
G.Z. is supported by a grant of NSFC and Qiushi Funding.}
\date{}

\maketitle
 \abstract
 We prove that the orbifold quantum ring is preserved
under  singular symplectic flops. Hence we verify Ruan's conjecture
for this case.
\endabstract
\def \rroot{e^{\frac{2\pi i}{r}}}
\def \mba{\mathbf{a}}
\def \mbb{\mathbf{b}}
\def \mbg{\mathbf{g}}
\def \mbh{\mathbf{h}}
\def \mbk{\mathbf{k}}
\def \inv{^{-1}}
\def \mzd{\mathbb{D}}

\tableofcontents

\section{Introduction}\label{section_1}

One of deep discovery in Gromov-Witten theory is its intimate
relation with the birational geometry. A famous conjecture of Ruan
asserts  that any two $K$-equivalent manifolds have isomorphic
quantum cohomology rings (\cite{R1}) (see also \cite{Wang}).
Ruan's conjecture was proved by Li-Ruan for smooth algebraic
3-folds (\cite{LR}) almost ten years ago. Only recently, it was
generalized to simple flops and Mukai flops in arbitrary
dimensions by Lee-Lin-Wang (\cite{LLW}). In a slightly different
context, there has been a lot of activities regarding Ruan's
conjecture in the case of McKay correspondence.

On the other hand, it is well known that the appropriate category
to study the birational geometry is not smooth manifolds. Instead,
one should consider the singular manifolds  with terminal
singularities. In the complex dimension three, the terminal
singularities are the finite quotients of hypersurface
singularities and hence the deformation of them are orbifolds. It therefore
raises the important questions if Ruan's conjecture still holds
for the orbifolds where there are several very interesting classes
of flops. This is the main topic of the current article.

Li-Ruan's proof of the case of smooth 3-folds consists of two
steps. The first step is to interpret flops in the  symplectic
category, then, they use almost complex deformation to reduce the
problem to the simple flop; the second step is to calculate the
change of quantum cohomology under the simple flop. The
description of a smooth simple flop is closely related to the
conifold singularity
$$
W_1=\{(x,y,z,t)|xy-z^2+t^2=0 \}.
$$

In \cite{CLZZ}, we initiate a program to understand the flop
associated with the singularities
$$
W_r=\{(x,y,z,t)|xy-z^{2r}+t^2=0\}/\mu_r(a,-a,1,0).
$$
$W_r$ appears in the list of terminal singularities in \cite{K}.
The singularities without quotient are also studied in
\cite{Laufer} and \cite{BKL}.

The program is along the same framework of that in \cite{LR}. The
first step is to describe the flops with respect to $W_r$
symplectically. This is done in the previous paper(\cite{CLZZ}).
Our main theorem in this paper is
\begin{theorem}\label{thm_0.1}
Suppose that $Y^s$ is a symplectic 3-fold with orbifold
singularities of type $W_{r_1},\ldots,W_{r_n}$ and $Y^{sf}$ is its singular flop,
then
$$
QH_{CR}(Y^s)=QH_{CR}(Y^{sf}).
$$
\end{theorem}

Theorem \ref{thm_0.1} verifies Ruan's conjecture in this
particular case. We should mention that Ruan also proposed a
simplified version of the above conjecture in terms of Ruan
cohomology $RH_{CR}$ which has been established in \cite{CLZZ} as
well. Furthermore, our previous results enters the proof of this
general conjecture in a crucial way.

The technique of the proof is a combination of the degeneration
formula of orbifold Gromov-Witten invariants, the localization
techniques and dimension counting arguments. The theory of relative orbifold Gromov-Witten
invariants and its degeneration formula involves heavy duty
analysis on moduli spaces and will appear elsewhere (\cite{CLS}).

The paper is organized as following. We first describe the
relative orbifold GW-invariants and state the degeneration formula
(without proof)(\S\ref{section_2}). Then, we summary the result of
\cite{CLZZ} on the singular symplectic flops and Ruan cohomology
(\S\ref{sect_4}). The heart of the proof is a detail analysis of
relative orbifold GW-theory on local models (\S \ref{sect_5} and
\S\ref{section_5}). The main theorem is proved in \S
\ref{section_6}.

{\em Acknowledge. }We would like to thank  Yongbin Ruan for
suggesting the problem and
 for many valuable discussions. We also wish to thank
Qi Zhang  for many discussions.

\section{Relative orbifold Gromov-Witten theory and the degeneration formula}\label{section_2}

\subsection{The Chen-Ruan Orbifold Cohomologies}\label{sect_1.1}

Let $X$ be an orbifold. For $x\in X$, if its small  neighborhood
$U_x$ is given by a uniformization system $(\tilde U, G, \pi)$, we
say $G$ is the isotropy group of $x$ and denoted by $G_x$. Let
$$
\mc T= \left(\bigcup_{x\in X} G_x\right)/\sim.
$$
Here $\sim$ is certain equivalence relation. For each $(g)\in \mc
T$, it defines a twisted sector $X_{(g)}$. At the mean while, the
twisted sector is associated with a degree-shifting number
$\iota(g)$. The Chen-Ruan orbifold cohomology is defined to be
$$
H^\ast_{CR}(X)= H^\ast(X) \oplus \bigoplus_{(g)\in \mc T}
H^{\ast-2\iota(g)}(X_{(g)}).
$$
For details, readers are referred to \cite{CR1}.
\def \OM{\overline{\mc M}}
\subsection{Orbifold Gromov-Witten invariants}\label{section_2.2}
Let $\OM_{g,n,A}(X)$ be the moduli space of representable orbifold
morphism of genus $g,n$-marked points and $A\in H_2(X,\mathbb{Z})$ (cf. \cite{CR2},\cite{CR3}).
By specifying the monodromy
$$
\mbh =((h_1),\ldots,(h_n))
$$
at each marked points, we can decompose
$$
\OM_{g,n,A}(X)=\bigsqcup_{\mbh} \OM_{g,n,A}(X,\mbh).
$$
Let
$$
ev_i:\OM_{g,n,A}(X,\mbh)\to X_{(h_i)}, 1\leq i\leq n
$$
be the evaluation maps. The primary orbifold Gromov-Witten
invariants are defined as
$$
\langle\alpha_1,\ldots, \alpha_n\rangle_{g,A}^X
=\int^{virt}_{\OM_{g,n,A}(X,\mbh)} \prod_{i=1}^m
ev_i^\ast(\alpha_i),
$$
where $\alpha_i\in H^\ast(X_{(h_i)})$.

In particular, set
\begin{eqnarray*}
&&\langle\alpha_1,\alpha_2, \alpha_3\rangle_{CR}
=\langle\alpha_1,\alpha_2, \alpha_3\rangle_{0,0},\\
&&\langle\alpha_1,\alpha_2, \alpha_3\rangle=
\langle\alpha_1,\alpha_2, \alpha_3\rangle_{CR}+
\sum_{A\not=0}\langle\alpha_1,\alpha_2, \alpha_3\rangle_{0,A}.
\end{eqnarray*}

\subsection{Ring structures on $H^\ast_{CR}(X)$. }\label{section_2.3}
Let $V$ be a vector space over $R$. Let
$$
h: V\otimes V\to R
$$
be a non-degenerate pairing and
$$
A: V\otimes V\otimes V\to R
$$
be a triple form. Then it is well known that one can define a
product $\ast$ on $V$ by
$$
h(u\ast v,w)=A(u,v,w).
$$
Different $A$'s give   different products.

\begin{remark}\label{remark_2.1}
Suppose we have $(V,h,A)$ and $(V',h',A')$. A map $\phi:V\to V'$
induces an isomorphism (with respect to the product) if $\phi$ is
a group isomorphism and
$$
\phi^\ast h'=h, \;\;\; \mbox{and} \;\;\; \phi^\ast A'=A.
$$
\end{remark}

Now let $V=H^\ast_{CR}(X)$ and $h$ be the Poincare pairing on $V$.
If
$$
A(\alpha_1,\alpha_2,\alpha_3)
=\langle\alpha_1,\alpha_2,\alpha_3\rangle_{CR},
$$
it defines the Chen-Ruan product. If
$$
A(\alpha_1,\alpha_2,\alpha_3)
=\langle\alpha_1,\alpha_2,\alpha_3\rangle,
$$
it defines the Chen-Ruan quantum product. We denote the ring to be
$QH^\ast_{CR}(X)$.

\subsection{Moduli spaces of relative stable maps for orbifold pairs}\label{section_2.4}
For the relative stable maps for the smooth case , there are two equivariant versions.
One is on the symplectic manifolds with respect to cylinder ends,
each of which admits a Hamiltonian $S^1$ action (\cite{LR}), the
other is on the closed symplectic manifolds with respect to
divisors(\cite{LR},\cite{Li}). This is also true for orbifolds. We adapt the
second version here.

Let $X$ be a symplectic orbifold with disjoint divisors
$$\{Z_1,\ldots,Z_k\}.$$
For simplicity, we assume $k=1$ and $Z=Z_1$.

By a relative stable map in $(X,Z)$, we mean a stable map  $f\in \mc
M_{g,n,A}(X)$ with additional data that  record how it intersects
with $Z$. Be precisely, suppose
$$
f: (\Sigma_g, \mf z)\to X.
$$
The additional data  are collected in order:
\begin{itemize}
\item Set
$$
\mbx=f\inv(Z)=\{x_1,\ldots, x_k\}.
$$
We call $x_i$ the {\em relative marked points}. The rest of marked
points are denoted by
$$
\mf p=\{p_1,\ldots,p_m\},
$$
i.e, $\mf z=\mbx\cup\mf p$;
\item Let
$$
\mbg=((g_1),\ldots, (g_k))
$$
denote the monodromy of $f$ (with respect to $Z$) at each point in
$\mbx$. The rest are denoted by
$$
\mbh=((h_1),\ldots, (h_m))
$$
which are the monodromy of $f$ (with respect to $X$) at each point
in $\mf p$.
\item the multiplicity of the tangency $\ell_j$ of $f$
with $Z$ at $z_j=f(x_j)$ is defined by the following: Locally, the
neighborhood of $z_j$ is given by
$$
(\tilde V\times \cplane\to \tilde V)/G_{z_j},
$$
where $\tilde V/G_{z_j}\subset Z$. Suppose the lift of $f$ is
\begin{eqnarray*}
&&\tilde f: \tilde \mzd\to \tilde V\times \cplane\\
&&\tilde f(t)=(v(t), u(t))
\end{eqnarray*}
Suppose the multiplicity of $u$ is $\alpha$ and  $g_j\in G_{z_j}$
acts on the fiber over $z_j$ with multiplicity $c$. Then the
multiplicity  is set to be
$$
\ell_j=\frac{\alpha\cdot c}{|g_j|}.
$$
We say $f$ maps $x_j$ to $Z_{g_j}$ at $\ell_j z_j$. Set
$$
\mf l =(\ell_1,\ldots,\ell_k).
$$
We may write
$$
f\inv(Z)=\mf l\cdot\mbx= \sum_{j=1}^k \ell_jx_j.
$$
\end{itemize}
As a relative object, we say $f$ is in the moduli space of relative
map
$$
\mc M_{g,n,A}(X,Z,\mbh,\mbg,\mf l).
$$
We denote the map by
$$
f: (\Sigma, \mf p,\mf l\cdot\mbx)\to (X,Z).
$$

We now describe the compactification of this moduli space. The
construction is similar to the smooth case(\cite{LR}).

The target space of a stable relative map is no longer $X$.
Instead, it is  extended in the following sense: let $L\to Z$ be
the normal bundle of $Z$ in $X$ and
$$
PZ=\mathbb{P}(L\oplus\cplane)
$$
 be its projectification, then
given an integer $b\geq 0$, we  have an extended target space
$$
X^\sharp_{b}:= X\cup \bigcup_{1\leq \alpha\leq b} PZ^\alpha.
$$
Here $PZ^\alpha$ denotes the $\alpha$-th copy of $PZ$. Let
$Z^\alpha_{0}$ be the 0-section and $Z^\alpha_{\infty}$ be the
$\infty$-section of $PZ^\alpha$. $X$ is called the root component
of $X^\sharp_b$. $Z_{0}^{b}$ is called the divisor
 of $X^\sharp_b$ and is (again) denoted by $Z$.

\begin{defn}\label{defn_2.2}
A relative map in $X^\sharp_{b}$ consists of following data: on
each component, there is a relative map: on the root component,
the map is denoted by
$$f^0: (\Sigma_0, \mf p^0, \mf l^0\cdot\mbx^0)\to (X,Z);$$
and on each component $PZ^\alpha$, the map is denoted by
$$
f^\alpha: (\Sigma^\alpha,\mf p^\alpha, \mf
l^\alpha\cdot\mbx^\alpha \cup\bar{\mf
l}^\alpha\cdot\bar{\mbx}^\alpha) \to (PZ^\alpha,Z_{0}^\alpha\cup
Z_\infty^\alpha).
$$
Here $\mf l^\alpha\cdot\mbx^\alpha=f\inv(Z_0^\alpha)$
 and  $\bar{\mf l}^\alpha\cdot\bar{\mbx}^\alpha=f\inv(Z_\infty^\alpha).$
Moreover, we require
 $f^\alpha$ at $Z_{0}^\alpha$ matches $f^{\alpha+1}$ at $Z_{\infty}^{\alpha+1}$. (see Remark \ref{remark_2.2}.)

We denote such a map by
$$
\mf f=(f^0,f^1,\ldots,f^b).
$$
Set $ \mbx=\mbx^{b}$ and
\begin{eqnarray*}
&&g_j=g_{x^b_j},\;\;\; \mbg=(g_1,\ldots,g_{|\mbx|})
\\
&&\ell_j=\ell_j^b, \;\;\; \mf l =(\ell_1,\ldots,\ell_{|\mbx|}).
\end{eqnarray*}
We say that $\mf f$ maps $x_{j}$ to the divisor $Z$ of
$X^\sharp_b$ at $\ell_jz_{j}\in Z_{(g_j)}$. Similarly, $\mbh$
records the twisted sector for
$$
\mf p=\mf p^0\cup\bigcup_{\alpha}\mf p^\alpha.
$$
The homology class $A$ in $X$ represented by $\mf f$ can be
defined properly. Collect the data
$$
\Gamma=(g,A,  \mbh,\mbg,\mf l),\;\;\; \mc T=( \mbg,\mf l).
$$
We say that $f$ is a relative orbifold map in $X^\sharp_{b}$ of
type  $(\Gamma,\mc T)$.  Denote the moduli space by
$\tilde\M_{\Gamma,\mc T}(X^\sharp_b,Z)$.
\end{defn}
\begin{remark}\label{remark_2.2}
Let $f^\alpha$ and $f^{\alpha+1}$ be as in the definition. Suppose
that
\begin{itemize}
\item $f^\alpha$ maps $x^\alpha_{j}\in \mbx^\alpha$
 to $(Z_{0}^\alpha)_{(g^\alpha_{j})}$
at $\ell^\alpha_kz^\alpha_{j}$; \item  $f^{\alpha+1}$
 maps $\bar{x}^{\alpha+1}_{i}\in \bar{\mbx}^{\alpha+1}$ to
 $(Z_{\infty}^{\alpha+1})_{(\bar{g}^{\alpha+1}_i)}$ at
$\bar{\ell}^{\alpha+1}_i \bar{z}^{\alpha+1}_{i}$,
\end{itemize}
 then by saying that {\em $f^\alpha$ at $Z_{0}^\alpha$
 matches $f^{\alpha+1}$ at $Z_{\infty}^{\alpha+1}$}
 we mean that
$$
\ell^\alpha_{i}=\bar{\ell}^{\alpha+1}_{i}, \;\;\;
z^\alpha_{i}=\bar{z}^{\alpha+1}_{i}, \;\;\; \mbox{and} \;\;\;
g^\alpha_{i}=\bar{g}^{\alpha+1}_i.
$$
\end{remark}
Note that there is a $\cplane^\ast $ action on $PZ^\alpha$. Let
$T$ be the product of these $b$ copies of $\cplane^\ast$. Then $T$
acts on $\tilde\M_{\Gamma,\mc T}(X^\sharp_b,Z)$. Define
$$
\M_{\Gamma,\mc T}(X^\sharp_b,Z)=\tilde\M_{\Gamma,\mc
T}(X^\sharp_b,Z)/T.
$$
It is standard to show that
\begin{prop}\label{prop_new_2.1}
There exists a large integer $B$ which depends on topological data $(\Gamma,\mc T)$
such that $\M_{\Gamma,\mc T}(X^\sharp_b,Z)$ is empty when $b\geq B$.
\end{prop}
Hence,
\begin{defn}\label{defn_2.3}
The compactified moduli space  is
$$
\om_{\Gamma,\mc T}(X,Z)=\bigcup_{b\in\integer^{\geq
0}}\M_{\Gamma,\mc T}(X^\sharp_b,Z).
$$
\end{defn}
The following technique theorem is proved in \cite{CLS}
\begin{theorem}\label{thm_2.1}
$\om_{\Gamma,\mc T}(X,Z)$ is a smooth compact virtual orbifold without
boundary with virtual dimension
$$
2c_1(A)+ 2(\dim_\cplane X-3)(1-g) + 2\sum_{i=1}^m(1-\iota(h_i)) +
2\sum_{j=1}^n(1-\iota(g_j)-[\ell_j]),
$$
where $[\ell_j]$ is the largest integer that is less or equal  to
$\ell_j$.
\end{theorem}

\subsection{Relative orbifold Gromov-Witten invariants.}\label{section_2.5}
Let $\om_{\Gamma,\mc T}(X,Z)$ be the moduli space given above.
There are evaluation maps
$$
ev_i^X:\om_{\Gamma,\mc T}(X,Z)\to X_{(h_i)}, \;\;\;
ev_i^X(f)=f(p_i), 1\leq i\leq m;
$$
and
$$
ev_j^Z:\om_{\Gamma,\mc T}(X,Z)\to Z_{(g_j)}, \;\;\;
ev_j^Z(f)=f(x_j), 1\leq j\leq k.
$$
Then for
$$\alpha_i\in H^\ast(X_{(h_i)}),1\leq i\leq m,\;\;\;
\beta_j\in H^\ast(Z_{(g_j)}),1\leq j\leq k
$$
the relative invariant is defined as
\begin{eqnarray*}
&&\langle \alpha_1,\ldots, \alpha_m|\beta_1,\ldots,\beta_k, \mc T\rangle_{\Gamma}^{(X,Z)}\\
&&\;\;\;\; =\frac{1}{|Aut(\mc T)|}\int_{\om_{\Gamma,\mc
T}(X,Z)}^{vir}\prod_{i=1}^m (ev_i^X)^\ast\alpha_i\prod_{j=1}^k
(ev_j^Z)^\ast\beta_j.
\end{eqnarray*}
In this paper, we usually set
$$
\mba=(\alpha_1,\ldots, \alpha_m),\;\;\;
\mbb=(\beta_1,\ldots,\beta_k),
$$
then the invariant is denoted by $\langle \mba|\mbb, \mc
T\rangle_{\Gamma}^{(X,Z)}$.

Moreover, if $\Gamma=\coprod_{\gamma}\Gamma^\gamma $, the relative
invariants (with disconnected domain curves) is defined to be  the
product of each connected component
$$
\langle \mba|\mbb, \mc T\rangle_{\Gamma}^{\bullet(X,Z)}
=\prod_\gamma\langle \mba|\mbb, \mc
T\rangle_{\Gamma^\gamma}^{(X,Z)}.
$$

\subsection{The degeneration formula}\label{section_2.6}
The symplectic cutting also holds for  orbifolds. Let $X$ be a
symplectic orbifold. Suppose that there is a local $S^1$
Hamiltonian action on $U\subset X$. We assume that
$$
U\cong Y\times (-1,1)
$$
and the projection onto the second factor
$$
\pi_2: U\to (-1,1)
$$
gives the Hamiltonian function. $Y\times\{0\}$ splits $X$ into two
orbifolds with boundary $Y$, denoted by $X^\pm$. Then the routine
symplectic cutting gives the degeneration
$$
\pi: X\to \bar X^+\cup_Z\bar X^-.
$$
Topologically, $\bar X^\pm$ is obtained by collapsing the
$S^1$-orbits of the boundaries of $X^\pm$.

There are maps
$$
\pi_\ast: H_2(X)\to H_2(X^+\cup_{Z} X^-), \;\;\; \pi^\ast:
H^\ast(X^+\cup_{Z} X^-)\to H^\ast(X).
$$
For $A\in H_2(M)$ we set $[A]\subset H_2(X)$ to be
$\pi_\ast\inv(\pi_\ast(A))$ and denote $\pi_\ast(A)$ by
$(A^+,A^-)$. On the other hand, for $\alpha^\pm\in H^\ast(X^\pm)$
with $\alpha^+|_Z=\alpha^-|_Z$, it defines a class on
$H^\ast(X^+\cup_{Z} X^-)$ which is denoted by
$(\alpha^+,\alpha^-)$. Let $\alpha=\pi^\ast(\alpha^+,\alpha^-)$.

\begin{theorem}\label{thm_3.1}
Suppose $\pi: X\to X^+\cup_{Z} X^-$ is the degeneration.
 Then
\begin{equation}\label{equation_2.1}
\langle \mba\rangle_{\Gamma}^X
=\sum_{I}\sum_{\eta=(\Gamma^+,\Gamma^-,I_\rho)} C_\eta\langle
\mba^+|\mbb^I,\mc T\rangle^{\bullet(X^+,Z)}_{\Gamma^+} \langle
\mba^-|\mbb_I,\mc T\rangle^{\bullet(X^-,Z)}_{\Gamma^-}.
\end{equation}
\end{theorem}
Notations in the formula are explained in order.
 $\Gamma$ is a data for Gromov-Witten invariants, it includes
$(g,[A])$; $(\Gamma^+,\Gamma^-,I_\rho)$ is an admissible triple
which consists of (possible disconnected) topological types
$\Gamma^\pm$ with the same relative data $\mc T$ under the
identification $I_\rho$ and they glue back to $\Gamma$. (For
instance,  one may refer to \cite{HLR} who interpret $\Gamma$'s as
graphs and then the gluing has an obvious geometric meaning); the
relative classes $\beta^i\in \mbb^I$ runs over a basis of
$Z_{(g_i)}$ and at the mean while $\beta_i$ runs over the dual
basis; finally
$$
C_\eta=|Aut(\mc T)|\prod_{i=1}^k\ell_i
$$
for $\mc T=( \mbg,\mf l)$.

\section{Singular symplectic flops}\label{sect_4}

\subsection{Local models and local flops}\label{sect_4.1}
Locally, we are concern those resolutions of
$$
\tilde W_r=\{(x,y,z,t)|xy-z^{2r}+t^2=0 \}
$$
and their quotients. $\tilde W_r-\{0\}$ inherits a symplectic form
$\tilde \omega_{r}^\circ$ from $\cplane^4$.

 By blow-ups, we have two small
resolutions of $\tilde W_r$:
\begin{eqnarray*}
\tilde{W}^s_r &=&\{((x,y,z,t),[p,q])\in \cplane^4\times \pone^1 \\
     && |xy-z^{2r}+t^2=0,\;\;\frac{p}{q}=\frac{x}{z^r-t}=\frac{z^r+t}{y}
     \}\\
\tilde{W}^{sf}_r &=&\{((x,y,z,t),[p,q])\in \cplane^4\times \pone^1 \\
     && |xy-z^{2r}+t^2=0,\;\;\frac{p}{q}=\frac{x}{z^r+t}=\frac{z^r-t}{y}
     \}.
\end{eqnarray*}

Let
$$
\tilde \pi_r^s:\tilde W^s_r\to \tilde W_r, \;\;\; \tilde
\pi_r^{sf}:\tilde W^{sf}_r\to \tilde W_r
$$
be the projections. The exceptional curves
$(\tilde\pi_r^s)^{-1}(0)$ and $(\tilde\pi_r^{sf})^{-1}(0)$ are
denoted by $\tilde\Gamma^s_r$ and $\tilde\Gamma^{sf}_r$
respectively. Both of them are isomorphic to $\pone^1$.

Let
$$
\mu_r=\langle\xi\rangle, \xi=\rroot
$$
be the cyclic group of $r$-th roots of 1. We denote its action on
$\cplane^4$ by $\mu_r(a,b,c,d)$ if the action is given by
$$
\xi\cdot(x,y,z,t)=(\xi^ax,\xi^by,\xi^cz,\xi^dt).
$$
Then $\mu_r(a,-a,1,0)$ acts on $\tilde W_r$, and naturally
extending to its small resolutions. Set
$$
W_r= \tilde W_r/\mu_r,\;\;\; W^s_r= \tilde W^s_r/\mu_r,\;\;\;
W^{sf}_r= \tilde W^{sf}_r/\mu_r.
$$
Similarly,
$$
\Gamma^s_r=\tilde{\Gamma}^s_r/\mu_r
\;\;\;\Gamma^{sf}_r=\tilde{\Gamma}^{sf}_r/\mu_r.
$$
We call that $W^s$ and $W^{sf}$  are the {\em small resolutions}
of $W_r$. We say that $W^{sf}$ is the flop of $W^s$ and vice
versa. They are both orbifolds with singular points on $\Gamma^s$
and $\Gamma^{sf}$. Note that the symplectic form
$\tilde\omega_r^\circ$ reduces to a symplectic form
$\omega_r^\circ$ on $W_r$.

It is known that
\begin{prop}\label{prop_4.1}
For $r\geq 2$, the normal bundle of $\tilde\Gamma^s_r$ {\em
($\tilde \Gamma^{sf}_r$)} in $\tilde W^s_r$ {\em ($\tilde
W^{sf}_r$)} is $\mc O\oplus \mc O(-2)$.
\end{prop}
{\bf Proof. } We take $\tilde W_r^s$ as an example. For the set
$$
\Lambda_p=\{q\not=0\},
$$
set $u=p/q$. Then $(u,z,y)$ gives a coordinate chart for
$\Lambda_p$. Similarly, for the set
$$
\Lambda_q=\{p\not=0\},
$$
set $v=q/p$. Then $(v,z,x)$ gives a coordinate chart for
$\Lambda_q$. The transition map is given by
\begin{equation}\label{eqn_4.1}
\left\{
\begin{array}{l}
v= u^{-1};\\
z=z;\\
x=-u^2y + 2uz^r.
\end{array}
\right.
\end{equation}
By linearize this equation, it is easy to get the conclusion.
q.e.d.
\begin{corollary}\label{cor_4.1}
For $r\geq 2$, the normal bundle of $\Gamma^s_r$ {\em
($\Gamma^{sf}_r$)} in $W^s_r$ {\em ($W^{sf}_r$)} is $(\mc O\oplus
\mc O(-2))/\mu_r$.
\end{corollary}

On $\tilde \Gamma^s_r$ ($\tilde \Gamma^{sf}_r$), there are two
special points. In term of $[p,q]$ coordinates, they are
$$
0=[0,1];\;\; \infty=[1,0].
$$
We denote them by $\mfp$ and $\mfq$ ($\mfP$ and $\mfQ$)
respectively. After taking quotients, they become singular points.
By the proof of Proposition \ref{prop_4.1},  the uniformization
system of $\mfp$ is
$$
\{ (p,x,y,z,t)|x=t=0\}
$$
with $\mu_r$ action given by
$$
\xi(p,y,z)=(\xi^ap,\xi^{-a}y,\xi z).
$$
At $\mfp$, for each given $\xi^k=\exp(2\pi ik/r), 1\leq k\leq r$,
there is a corresponding twisted sector(\cite{CR1}). As a set, it
is same as $\mfp$. We denote this twisted sector by $[\mfp]_{k}$.
For each twisted sector, a degree shifting number is assigned. We
conclude that
\begin{lemma}\label{lemma_4.1}
For $\xi^k=\exp(2\pi ik/r), 1\leq k\leq r$, the degree shifting
$$
\iota([\mfp]_{k})= 1+ \frac{k}{r}.
$$
\end{lemma}
\n{\bf Proof. }This follows directly from the definition of degree
shifting. q.e.d.

\v Similar results hold for the singular point $\mfq$. Hence we
also have twisted sector $[\mfq]_k$ and
$$
\iota([\mfq]_{k})= 1+ \frac{k}{r}.
$$
Similarly, on $W^{sf}$, there are twisted sectors $[\mfP]_k,
[\mfQ]_k$ and
$$
\iota([\mfP]_k)=\iota([\mfQ]_k)=1+\frac{k}{r}.
$$

\subsection{Torus action.}\label{section_3.2}
We introduce a $T^2$-action on $\tilde W_r$:
$$
(t_1,t_2)(x,y,z,t)= (t_1t_2^{r}x,t_1^{-1}t_2^{r}y, t_2z, t_2^{r}t).
$$
For an action $t_1^{a} t_2^{b}\cdot$, we write the weight of
action by $a\lambda+bu$. This action naturally extends to the
actions on all models generated from $\tilde W_r$, such as
$W_r,W^{s}_r$ and $W_r^{sf}$.

It then induces an action on $\tilde \Gamma^s_r$
($\tilde{\Gamma}^{sf}_r$):
$$
(t_1,t_2)[p,q]= [t_1p,q].
$$
Recall that the normal bundle of $\tilde\Gamma^s_r$ in $\tilde
W^s_r$ is $\mc O\oplus \mc O(-2)$.
\begin{lemma}\label{lemma_5.2}
The action weights at $\mc O_p$ and $\mc O_q$ are $u$. The action
weights at $\mc O_p(-2)$ and $\mc O_q(-2)$ are $-\lambda +ru$ and
$\lambda + ru$.
\end{lemma}
{\bf Proof. }This follows directly from the model given by
\S\ref{sect_4.1}. q.e.d.

\v\n
 it is easy to verify that
$\mfp,\mfq$ ($\mfP,\mfQ$) are fixed points of the action.

On the other hand, there are four special lines connecting to
these points that are {\em invariant} with respect to the action.
Let us look at $\tilde W^s_r$. For the point $\mfp$, two lines are
in $\Lambda_p$ and are  given by
\begin{eqnarray*}
\tilde L^s_{p,y}& =&\{x=z=t=0,u=0\},\\
\tilde L^s_{p,z}&=&\{x=y=0, z^r+t=0,u=0\}.
\end{eqnarray*}
To the point $\mfq$, two lines are in $\Lambda_q$ and are  given
by
\begin{eqnarray*}
\tilde L^s_{q,x}& =&\{y=z=t=0,v=0\},\\
\tilde L^s_{q,z}&=&\{x=y=0, z^r-t=0,v=0\}.
\end{eqnarray*}
Similarly, for $\tilde W^{sf}_r$  we have
\begin{eqnarray*}
\tilde L^{sf}_{p,y}& =&\{x=z=t=0,u=0\},\\
\tilde L^{sf}_{p,z}&=&\{x=y=0, z^r-t=0,u=0\},\\
\tilde L^{sf}_{q,x}& =&\{y=z=t=0,v=0\},\\
\tilde L^{sf}_{q,z}&=&\{x=y=0, z^r+t=0,v=0\}.
\end{eqnarray*}
Correspondingly, these lines in $W^s_r$ and $W^{sf}_r$ are denoted
by the same notations without tildes.

\begin{remark}\label{remark_3.1}
Note that the defining equations for the pairs $L^s_{q,x}$ and
$L^{sf}_{q,x}$, $L^s_{p,y}$ and $L^{sf}_{p,y}$ are same.
\end{remark}

\subsection{Symplectic orbi-conifolds and singular symplectic flops }\label{sect_4.2}
An orbi-conifold (\cite{CLZZ}) is a topological space $\mc Z$ with
a set of (singular) points
$$
P=\{p_1,\ldots,p_k\}
$$
such that $\mc Z-P$ is an orbifold and for each $p_i\in P$ there
exists a neighborhood $U_i$ that is isomorphic to $W_{r_i}$ for
some integer $r_i\geq 1$. By a symplectic structure on $\mc Z$ we
mean a symplectic form $\omega$ on $\mc Z-P$ and it is
$\omega^\circ_{r_i}$ in $U_i$. We call $\mc Z$ a symplectic
orbi-conifold.
 There exists $2^k$ resolutions of $\mc Z$. Let $Y^s$ be such a resolution, its flop
is defined to be the one
 that is obtained by flops each local model of $Y^s$. We denote it by $Y^{sf}$.
In \cite{CLZZ} we prove that
\begin{theorem}\label{thm_4.2}
 $Y^s$ is a symplectic orbifold if and only if $Y^{sf}$ is.
\end{theorem}
So $Y^{sf}$ is called the (singular) symplectic flop of $Y^s$ and
vice versa.

Now for simplicity, we assume that $\mc Z$ contains only one
singular point $p$ and is smooth away from $p$. Suppose $Y^s$ and
$Y^{sf}$ are two resolutions that are flops of each other and,
locally, $Y^s$ contains $W^s_r$ and $Y^{sf}$ contains $W^{sf}_r$.
Then
\begin{eqnarray*}
&&H^\ast_{CR}(Y^s)= H^\ast(Y^s)\oplus\bigoplus_{k=1}^r
\cplane[\mfp]_k \oplus\bigoplus_{k=1}^r\cplane[\mfq]_k;
\\
&& H^\ast_{CR}(Y^{sf})= H^\ast(Y^{sf})\oplus\bigoplus_{k=1}^r
\cplane[\mfP]_k \oplus\bigoplus_{k=1}^r\cplane[\mfQ]_k.
\end{eqnarray*}
\begin{lemma}\label{lemma_4.2}
There are natural isomorphisms
$$
\psi_k:H^k(Y^s)\to H^{k}(Y^{sf}).
$$
\end{lemma}
{\bf Proof. }We know that
$$
Y^s-\Gamma^s=Y^{sf}-\Gamma^{sf}.
$$
We also have  the exact sequence
$$
\cdots \to H^k(Y,Y\setminus \Gamma)\to H^k(Y)\to
H^{k}(Y\setminus\Gamma)\to H^{k+1}(Y,Y\setminus \Gamma)\to\cdots
$$
and
$$
H^{k}(Y,Y\setminus\Gamma)\cong H^k_c(\mc O\oplus \mc O(-2)) \cong
H^{k-4}(\pone^1).
$$
$Y$ is either $Y^s$ or $Y^{sf}$ and $\Gamma$ is the exceptional
curve in $Y$.

Suppose we have $\omega^s\in H^k(Y^s)$. Suppose $X^s_{\omega^s}$
is its Poincare dual. If $k>2$, we may require that
$X^s_{\omega^s}\cap \Gamma^s=\emptyset$. Hence $X^s_{\omega^s}$ is
in
$$
Y^s\setminus \Gamma^s=Y^{sf}\setminus\Gamma^{sf}.
$$
Using this, we get a class $\omega^{sf}\in H^k(Y^{sf})$. Set
$\psi_k(\omega^s)=\omega^{sf}$.

If $k\leq 2$, since
$$
H^m_{comp}(\mc O\oplus \mc O(-2))=0, m\leq 3,
$$
we have
$$
H^k(Y^s)\cong H^k(Y^s\setminus \Gamma^s)\cong H^k(Y^{sf}\setminus
\Gamma^{sf}) \cong H^k(Y^{sf}).
$$
The isomorphism gives $\psi_k$. q.e.d.

\v\n On the other hand, we set
$$
\psi_o([\mfp]_k)=[\mfP]_k, \;\;\;\psi_o([\mfq]_k)=[\mfQ]_k,
$$
Totally, we combine $\psi_k$ and $\psi_o$ to get a map
\begin{equation}\label{eqn_4.2}
\Psi^\ast: H^\ast_{CR}(Y^s)\to H^\ast_{CR}(Y^{sf}).
\end{equation}
It can be shown that
\begin{prop}\label{prop_4.2}
$\Psi^\ast$ preserves the Poincare pairing.
\end{prop}
Without considering the extra classes from twisted sectors, the proof is standard.
When the cohomology classes from twisted sectors  are involved,
it is proved in \cite{CLZZ}.

On the other hand, there is a natural isomorphism
$$
\Psi_\ast: H_2(Y^s)\to H_2(Y^{sf})
$$
with $\Psi_\ast([\Gamma^s_r])=-[\Gamma^{sf}_r]$.

Now suppose  that we do the symplectic cutting on $Y^s$ and
$Y^{sf}$ at $W^s_r$ and $W^{sf}_r$ respectively. Then
\begin{eqnarray*}
&&\pi_s: Y^s\xrightarrow{degenerate} Y^-\cup_Z M^s_r;\\
&&\pi_{sf}:Y^{sf}\xrightarrow{degenerate} Y^-\cup_Z M^{sf}_r.
\end{eqnarray*}
It is clear that $M^s_r$ and $M^{sf}_r$ are flops of each other.
Then similarly, we have a  map
$$
\Psi^\ast_r: H^\ast_{orb}(M^s_r)\to H^\ast_{orb}(M^{sf}_r).
$$
It is easy to see that the diagram
\begin{equation}\label{eqn_4.3}
\begin{CD}
H^{\ast}_{CR}(Y^-\cup_Z M^s_r)   @>(id, \sigma^\ast)>>  H^{\ast}_{CR}(Y^-\cup_Z M^{sf}_r)\\
@V{\pi_s^\ast} VV      @VV{\pi_{sf}^\ast}V \\
H^\ast_{CR}(Y^s)   @>\Sigma^\ast>> H^\ast_{CR}(Y^{sf})
\end{CD}
\end{equation}
commutes.

\subsection{Ruan cohomology rings}\label{section_3.4}
As explained in \S\ref{section_2.3}, the cohomology ring structure is defined via
a triple form $A$. In the current situation, we can define  a ring
structure on $Y^s$ (and $Y^{sf}$) that plays a role between
Chen-Ruan (classical) ring structure and Chen-Ruan quantum ring
structure. The triple forms on $Y^s$ and $Y^{sf}$ are given by
\begin{equation}
\langle\alpha_1,\alpha_2, \alpha_3\rangle_R=
\langle\alpha_1,\alpha_2, \alpha_3\rangle_{CR}+
\sum_{A=d[\Gamma_r^s],d>0}\langle\alpha_1,\alpha_2,
\alpha_3\rangle_{0,A}q_s^d,
\end{equation}
\begin{equation}
\langle\alpha_1,\alpha_2, \alpha_3\rangle_R=
\langle\alpha_1,\alpha_2, \alpha_3\rangle_{CR}+
\sum_{A=d[\Gamma_r^{sf}],d>0}\langle\alpha_1,\alpha_2,
\alpha_3\rangle_{0,A}q_{sf}^d
\end{equation}
respectively. Here $q_s$ and $q_{sf}$ are formal variables that
represent classes $[\Gamma_r^s]$ and $[\Gamma_{r}^{sf}]$. They
define  Ruan rings  $RH(Y^s)$ and $RH(Y^{sf})$. In \cite{CLZZ}, we
already proved that
\begin{theorem}\label{theorem_3.1}
$\Psi^\ast$ gives the isomorphism $RH(Y^s)\cong RH(Y^{sf})$.
\end{theorem}

\section{Relative Gromov-Witten theory on $M^s_r$ and $M^{sf}_r$}\label{sect_5}

\subsection{Local models $M^s_r$ and $M^{sf}_r$}\label{sect_5.1}
$M^s_r$ and $M^{sf}_r$ are obtained from $W^s_r$ and $W^{sf}_r$ by
cutting at infinity. We explain this precisely.

We introduce an $S^1$ action on $\cplane^4$
$$
\gamma(x,y,z,t)=(\gamma^r x,\gamma^r y,\gamma z,\gamma^r t).
$$
Using this action, we collapse $\tilde W_r$ at $\infty$. The
infinity divisor is identified as
$$
\tilde Z=\frac{\tilde W^s_r\cap S^7}{S^1}.
$$
By this way, we get an orbifold with singularity at 0, denoted by
$\tilde M_r$. By blowing-up $\tilde M_r$ at 0, we have $\tilde
M_r^s$ and $\tilde M^{sf}_r$. $\mu_r$-action  can then naturally
extend to $\tilde M_r$, $\tilde M^s_r$ and $\tilde M^{sf}_r$. By
taking quotients, we have $M_r$, $ M^s_r$ and $M^{sf}_r$. $M^s_r$
and $M^{sf}_r$ are the collapsing of $W^s_r$ and $W^{sf}_r$ at
infinity. Note that the $T^2$-action given in \S\ref{section_3.2} also acts on
these spaces.

\def \mbP{\mathbb{P}}
Let $\tilde\mbP:=\mathbb{P}(r,r,1,r,1)$ be the weighted projective
space. Then $\tilde M_r$ can be embedded in $\tilde\mbP$ and is
given by the equation
$$
xy-z^{2r}+t^2=0.
$$
The original $\tilde W_r$ is embedded in $\{w\not=0\}$ and $\tilde
Z$ is in $\{w=0\}$. $\mu_r$-action extends to $\tilde\mbP$ by
$$
\xi(x,y,z,t,w)=(\xi^ax,\xi^{-a}y,\xi z,t,w).
$$
Then $M_r$ is embedded in $\mbP:=\tilde\mbP/\mu_r$. Set $Z=\tilde
Z/\mu_r$. To understand the local behavior of $M^s_r$ and
$M^{sf}_r$ at $Z$, it is sufficient to use this model at
$\{w=0\}$.

\def \mfx{\mathfrak{x}}
\def \mfy{\mathfrak{y}}
\def \mfzp{\mathfrak{z}^+}
\def \mfzm{\mathfrak{z}^-}
We now study the singular points at $Z$. Combining the $S^1$ and
$\mu_r$ actions, we have
$$
(\gamma,\xi)(x,y,z,t,w) =(\gamma^r\xi^a x,
\gamma^r\xi^{-a}y,\gamma\xi z, \gamma^r t,\gamma w),
(\gamma,\xi)\in S^1\times \mathbb{Z}_r.
$$
\begin{lemma}\label{lemma_5.1}
There are four singular points
\begin{eqnarray*}
&&\mfx=[1,0,0,0,0];\\
&&\mfy=[0,1,0,0,0];\\
&&\mfzp=[0,0,1,1,0];\\
&&\mfzm=[0,0,1,-1,0]
\end{eqnarray*}
and a singular set
$$
S=\{xy+t^2=0, z=0\}
$$
on $Z$. Their stabilizers are $\mathbb{Z}_{r^2}, \mathbb{Z}_{r^2},
\mathbb{Z}_{r}$, $\mathbb{Z}_{r}$ and $\mathbb{Z}_{r}$
respectively.
\end{lemma}
{\bf Proof. }Take a point $(x,y,z,t,w)$. We find those points with
nontrivial stabilizers. \v\n Case 1, assume that $z\not=0$. Then
$\gamma\xi=1$. Therefore,
$$
(\gamma,\xi)(x,y,z,t,w) =(\gamma^r\xi^a x,
\gamma^r\xi^{-a}y,\gamma\xi z, \gamma^r t,0)=(\xi^a x,\xi^{-a}y,z,
t).
$$
In order to have nontrivial stabilizers, we must have $x=y=0$.
Therefore, only $\mfzp$ and $\mfzm$ survive. Their stabilizer are
both $\mathbb Z_r$. \v\n Case 2, assume that $z=0$. If $t\not=0$,
$\gamma$ should be a $r$-root. Then
$$
(\gamma,\xi)(x,y,z,t,w) =(\xi^a x, \xi^{-a}y,0, t,0)=(\xi^a
x,\xi^{-a}y,z, t).
$$
Hence, when $xy$ are not 0, the set $\{xy+t^2=0\}$ has  the
stabilizer $\mathbb{Z}_{r}$.

\v\n Case 3, assume that $z=t=0$. Then $xy=0$ by the equation.
Hence we can only have $\mfx$ and $\mfy$. Clearly, their
stabilizers are both $\mathbb{Z}_{r^2}$. q.e.d.

\v\n We now look at the local models at $\mfzp$ and $\mfzm$. The
coordinate chart at $\mfzp$ is given by $(x,y,w)$ and the action
is
\begin{equation}\label{eqn_new_1}
\xi(x,y,w)= (\xi^{-a} x,\xi^a y, \xi w), \xi\in \mathbb{Z}_r.
\end{equation}
The model at $\mfzm$ is same.

Now look at the local models at $\mfx$ and $\mfy$. At $\mfx$, the
local coordinate chart is given by $(z,t,w)$. The action is given
by
\begin{equation}\label{eqn_new_2}
\xi(z,t,w)= (\xi\eta z, \xi^r t,\xi w), \;\;\; \mbox{where}\;\;\;
\eta^a\xi^r=1, \xi\in \mathbb{Z}_{r^2}.
\end{equation}
Suppose
\begin{equation}\label{equation_4.2}
\eta= \exp{2\pi i \mu}, 0\leq \mu<1.
\end{equation}
Similarly, at $\mfy$, the local coordinate chart is given by
$(z,t,w)$. The action is given by
\begin{equation}\label{eqn_new_3}
\xi(z,t,w)= (\xi\eta z, \xi^r t,\xi w), \;\;\; \mbox{where}\;\;\;
\eta^{-a}\xi^r=1,\xi\in \mathbb{Z}_{r^2}.
\end{equation}

For points on $S$, the action on the normal direction is given by
\begin{equation}\label{eqn_new_4}
\xi(z,w)=(\xi z,\xi w), \xi\in \mathbb{Z}_r.
\end{equation}

Recall that we have four lines described in \S\ref{section_3.2}. They
are now being four lines in $M^s_r$ ($M^{sf}_r$). Take $M^s_r$ as
an example. We have
\begin{eqnarray*}
& L^s_{p,y}: \mfp\leftrightarrow \mfy;
& L^s_{q,x}: \mfq\leftrightarrow \mfx; \\
& L^s_{p,z}: \mfp\leftrightarrow \mfzm; & L^s_{q,z}:
\mfq\leftrightarrow \mfzp;.
\end{eqnarray*}
In the table, for each line we give the name of the curve and the
ends it connects.

For $M^{sf}_r$, we have
\begin{eqnarray*}
& L^{sf}_{p,y}: \mfP\leftrightarrow \mfy;
& L^{sf}_{q,x}: \mfQ\leftrightarrow \mfx; \\
& L^{sf}_{p,z}: \mfP\leftrightarrow \mfzp; & L^{sf}_{q,z}:
\mfQ\leftrightarrow \mfzm;.
\end{eqnarray*}

Regarding the $T^2$-action, we have following two lemmas. The
proof is straightforward, we leave it to readers.

 \begin{lemma}\label{lemma_5.3}
The fixed points on $M^s_r$ ($M^{sf}_r$) are $\mfp$, $\mfq$
($\mfP,\mfQ$) on $\Gamma^s_r$ ($\Gamma^{sf}_r$) and
$\mfx,\mfy,\mfzp,\mfzm$ on $Z$.
\end{lemma}

\begin{lemma}\label{lemma_5.4}
In $M^s_r$, the invariant curves with respect to the torus action
are $L^s_{p,y},L^s_{q,x},L^s_{p,z}$ and $L^s_{q,z}$.
\end{lemma}

\subsection{Relative Moduli spaces for the pair $(M^s_r,Z)$}\label{sect_5.3}

We explain the relative moduli spaces for the pair $(M^s_r,Z)$.
Similar explanations can be applied to $(M^{sf}_r,Z)$.

Let
$$
\Gamma=(g,A,\mbh,\mbg,\mf l),\;\;\; \mc T=(\mbg,\mf l)
$$
be as before. Let $\om_{\Gamma,\mc T}(M^s_r,Z)$ be the moduli
space. Recall that the virtual dimension of the moduli space is
$$
\dim=c_1(A) +\sum_{i=1}^m(1-\iota(h_i))
+\sum_{j=1}^k(1-[\ell_{j}]-\iota(g_{j})).
$$
Here we use the complex dimension.

First, we note that
$$
[c_1(M^s_r)\cdot Z]= (r+2)[c_1(L_Z)\cdot Z]=\sum_{x\in\mbx}(r+2)\ell_x.
$$
Therefore
$$
\dim=\sum_{i=1}^m(1-\iota(h_i))
+\sum_{j=1}^k((r+2)\ell_j+1-[\ell_{j}]-\iota(g_{j}))
$$
\def \mfk{\mathfrak}
Set
\begin{eqnarray*}
&&\mfk u_i= 1-\iota(h_i), 1\leq i\leq m,\\
&&\mfk v_j=(r+2)\ell_j+1-[\ell_{j}]-\iota(g_{j}), 1\leq j\leq k.
\end{eqnarray*}
Then
$$
\dim=\sum_{i=1}^m\mfk u_i+\sum_{j=1}^k\mfk v_j.
$$
For $\mfk u_i$ and $\mfk v_j$ we have following facts:
\begin{enumerate}
\item if $x_j$ is mapped to $\mfx$ with
 $g_j=\exp(2\pi i\frac{\alpha}{r^2} )$. Then
$$
\ell_{j}=[\ell_j]+ \frac{\alpha}{r^2}.
$$
So
$$
\mfk v_j= (r+2) (\frac{\alpha}{r^2}+[\ell_j])+1-([\ell_j] +
\frac{\alpha}{r^2}+ \{\frac{\alpha}{r^2}+\mu\}+
\{\frac{\alpha}{r}\}).
$$
Here $\{z\}:=z-[z]$. Note that this can be
$$
(r+1)[\ell_j]+n -\mu, n\geq 1.
$$
Here $\mu$ is defined in \eqref{equation_4.2}.
The degree shifting numbers are given by \eqref{eqn_new_2}.
\item If $x_j$ is mapped to
$\mfy$ with
 $g_j=\exp(2\pi i\frac{\alpha}{r^2} )$, it is same as the previous case.
\item If $x_j$ is mapped to $S$ with $g_j=\exp(2\pi
i\frac{\alpha}{r} )$, then
$$
\mfk v_j=(r+1)[\ell_j]+\alpha+1.
$$
The degree shifting numbers are given by \eqref{eqn_new_4}.
\item  If $x_j$ is mapped to $\mfzp$ with
 $g_j=\exp(2\pi i\frac{\alpha}{r} )$.
Then
$$
\ell_j=[\ell_j]+ \frac{\alpha}{r}.
$$
So
$$
\mfk v_j= (r+1)[l_x] +\alpha +\frac{\alpha}{r}.
$$
Here the degree shifting numbers are given by \eqref{eqn_new_1}.
\item If $x_j$ is mapped to $\mfzm$ with
 $g_j=\exp(2\pi i\frac{\alpha}{r} )$, then it is same as the previous case.

\item whenever $g_j=1$
$$
\mfk v_j=(r+1)\ell_j +1.
$$
\item  when  $h_i= \exp(2\pi i\frac{\alpha}{r})$, $ \mfk
u_i=-\frac{\alpha}{r}. $ \item  when $h_i=1$, $\mfk u_i=1$.
\end{enumerate}

\subsection{Admissible data $(\Gamma,\mc T,\mba)$}\label{sect_5.4}
Suppose we are computing the  relative Gromov-Witten invariant
\begin{equation}\label{equation_4.1}
\langle \mba|\mbb,\mc T\rangle_{\Gamma}^{(M^s,Z)}.
\end{equation}

Let $|\alpha|$ denote the degree of a form $\alpha$. Set
$$
N= \dim-\sum_{i=1}^m|\alpha_i|.
$$
On the other hand,  set
$$
N'= \sum_{j=1}^k \dim (Z_{g_j}).
$$
\begin{defn}\label{defn_5.1}
The data $(\Gamma,\mc T,\mba)$ is called admissible  if $N\leq
N'$.
\end{defn}
By the definition, we have
\begin{lemma}\label{LEMMA_4.1}
If $(\Gamma,\mc T,\mba)$ is not admissible, the invariant \eqref{equation_4.1} is 0.
\end{lemma}
In this paper, we may assume that:
\begin{assumption}\label{assump}
(i) $|\alpha_i|=0$ for all
$p_i$, (ii) $|\mba|\leq 3$.
\end{assumption}
Since
$$
N-N'=\sum_{i=1}^m\mfk u_i+\sum_{j=1}^k (\mfk v_j-\dim (Z_{g_j})),
$$
we make the following definition.
\begin{defn}\label{definition_4.1}
 we say that $\mfk u_i$ is the contribution of marked point $p_i$ to $N-N'$ and
 $\mfk v_j-\dim (Z_{g_j})$ is that of $x_j$.
\end{defn}

\begin{prop}\label{prop_5.2}
Suppose that Assumption \ref{assump} holds. If $(\Gamma,\mc T,\mba)$ is admissible, then  one of the following
cases holds:
\begin{enumerate}
\item $\mbx$ consists of only one smooth point, then $\mf p$
consists of three singular points $(p_1,p_2,p_3)$ such that
$$
\iota(h_{1})+\iota(h_{2})+\iota(h_{3})=5.
$$
For this case, $N=N'$. \item if $\mbx$ contains a point $x$ maps
to $\mfzm$ or $\mfzp$, then one of the following should hold:
\begin{eqnarray*}
&&|\mba|=2, \iota(h_1)+\iota(h_2)= 3+\frac{1}{r};\\
&&|\mba|=3, \iota(h_1)+\iota(h_2)+\iota(h_3)= 4+\frac{1}{r}.
\end{eqnarray*}
\item $\mbx$ consists of only singular points, the multiplicities
at $x\in \mbx$ are all less than 1. Furthermore, $x$ can not be
mapped to either $\mfzp$ or $\mfzm$.
\end{enumerate}
Moreover, $\ell_x\leq 1$.
\end{prop}
{\bf Proof. } First, we suppose that $x$ is smooth point. It
contributes $(r+1)l_x+1$ to $N$ and contributes $2$ to $N'$. Hence
its contribution to $N-N'$ is
$$
(r+1)l_x-1\geq r.
$$

If $x$ a singular point, its contribution to $N-N'$ is given by
the following list
\begin{itemize}
\item $x\to \mfx$ but not in $S$, the contribution is
$(r+1)[l_x]+n -\mu, n\geq 1$;
\item $x\to S$, the contribution is $ (r+1)[l_x]+\alpha+1, $ \item
$x\to \mathfrak{z}^{\pm}$, the contribution is $(r+1)[l_x]+\alpha +\alpha/r$.
\end{itemize}
Note that they are all positive. Hence we conclude that, if $\mbx$
contains a smooth point $x$, then only the following situation
survives: $|\mbx|=1$, $r=2, |\mba|=3$, and
$$
\iota(h_{1})+\iota(h_{2})+\iota(h_{3})=5.
$$
Furthermore, $N=N'$.

Now suppose that  $\mbx$ contains a point $x$ mapping to $\mfzm$
(or $\mfzp$), only the following situation survives: $\alpha=1$
and one of the following holds:
\begin{eqnarray*}
&&|\mba|=2, \iota(h_1)+\iota(h_2)= 3+\frac{1}{r};\\
&&|\mba|=3, \iota(h_1)+\iota(h_2)+\iota(h_3)= 4+\frac{1}{r}.
\end{eqnarray*}
The rest admissible data belong to the third case. q.e.d.

\section{Vanishing results on relative invariants}\label{section_5}

\subsection{Localization via the torus action}\label{sect_5.5}

The torus action $T^2$ on $M^s_r$ induces an action on the moduli
space. We now study the fix loci of the moduli space with respect
to the action. We use the notations in \S\ref{section_2.4} for
$X=M^s_r$.

A relative stable map
$$
\mf f=(f^0,f^1_1,\ldots,f^{b_1}_1,\ldots, f^k_1,\ldots, f^{b_k}_k)
$$
is invariant if and only if each $f$ is invariant.
Since we only consider the invariants for admissable data,
the invariant maps in such moduli spaces only have $f^0$.
In fact, the fact $\ell_x\leq 1$ in Proposition \ref{prop_5.2},
implies this.

$f^0$ is a stable map in $X$
whose components are invariant
maps (maybe constant map) and nodal points are mapped to fix
points, which are $\mfp,\mfq$ and
$\mfx,\mfy,\mfzp,\mfzm$.
The constant map should also map to these
points, while the nontrivial invariant curves should cover one of
those four lines or $\Gamma^s_r$. Set
$$
\mbox{FT}=\{\mfp,\mfq,\mfx,\mfy,\mfzp,\mfzm\},\;\;\;
\mbox{IC}=\{\Gamma^s_r,L^s_{p,y},L^s_{q,x},L^s_{p,z}, L^s_{q,z}\}
$$

As in the Gromov-Witten theory, we introduce graphs to describe
the components of fix loci. We now describe the graph $T$ for
$f^0$. Let $V_T$ and $E_T$ be the set of vertices and edges of
$T$.
\begin{itemize}
\item each vertex is assigned to a connected component of the
pre-image of  $\mbox{FT}$; on each vertex, the image point is
recorded; \item each edge is assigned to the component that is
non-constant map; the image with multiplicity is recorded; \item
on each flag, a twisted sector (or the group element of the
sector) is recorded.
\end{itemize}
Let $F_T$ be the fix loci that correspond to the graph $T$. Since
$T$ only describes $f^0$, $F_T$ may contain several components.
Let $\mathfrak{T}$ be the collection of graphs and $\mc F$ be the
collection of $F_T$.

We recall the virtual localization formula. Suppose that $\Omega$
is a form on  $\om_\Delta(X,Z)$ and $\Omega_{T_2}$ is its
equivariant extension if exists. Then
$$
I(\Omega)=\int_{\om_{\Delta}(X,Z )}^{vir}
\Omega=\sum_{T\in\mathfrak{
T}}\int_{F_T}\frac{\Omega_{T^2}}{e_{T^2}(N^{vir}_{F_T})}.
$$
Here, $N^{vir}_{F_T}$ is the virtual normal bundle of $F_T$ in the
virtual moduli space, $e_{T^2}$ is the $T^2$-equivariant Euler
class of the bundle. The right hand side is a function in
$(\lambda,u)$, which is rational in $\lambda$ and polynomial in
$u$. We denote each term in the summation as
$I_{F_T}^\Omega(\lambda,u)$ and the sum by $I^\Omega(\lambda,u)$.

\begin{lemma}\label{LEMMA_5.1}
Let $(\Gamma,\mc T,\mba)$ be an
admissible data in Proposition \ref{prop_5.2}.
Then the nontrivial relative invariants
$$
\langle \mba|\mbb\rangle_{\Gamma,\mc T}^{(M^s,Z)}
$$
can be computed via localization.
\end{lemma}
{\bf Proof. }It is sufficient to show that the forms in $\mba$ and
$\mbb$ have equivariant extensions. By our assumption, we always
take $\alpha\in \mba$ to be 1. It is already equivariant.

Now suppose that $\beta_j$ is assigned to $x_j$. If $Z_{(g_j)}$ is
a single point, $\beta_j$ are taken to be 1, which is equivariant.
If $Z_{g_j}=S$, $\beta_j$ is either a 0 or 2-form, both have
equivariant representatives. The last case is that $x$ is smooth.
For this case, since $N=N'$ (cf. case (1) in Proposition
\ref{prop_5.2}),  we must have $\deg(\beta)=4$ to get
nontrivial invariants. Since $\beta$ is of top degree, it has an
equivariant representative as well. q.e.d.

\subsection{Vanishing results on $I^\Omega_{F_T}(\lambda,u)$, (I)}\label{sect_5.6}
By localization, we have
$$
I(\Omega)=I^\Omega(\lambda,u)=\sum_{T\in
\mathfrak{T}}I^\Omega_{F_T}(\lambda,u).
$$
Since the left hand side is independent of $u$, we have
$$
I(\Omega)=\lim_{u\to 0}I^\Omega(\lambda,u)=\sum_{T\in
\mathfrak{T}}\lim_{u\to 0} I^\Omega_{F_T}(\lambda,u).
$$

\begin{theorem}\label{prop_5.3}
Suppose that  $|\mbx|>0$.
If $T$ contains an edge $e_0$ that
records a map cover  $\Gamma^s_r$, then
$$
\lim_{u\to 0}I_{F_T}(\lambda,u)=0.
$$
\end{theorem}
Edge $e_0$  records a map:
$$
f_0:S^2\to \Gamma^s_r.
$$
$f_0$ can be either a smooth or an orbifold map. Hence, we restate
the theorem as,
\begin{prop}\label{lemma_5.5}
If the map $f_0$ for $e_0$ is smooth, then
$$
\lim_{u\to 0}I_{F_T}(\lambda,u)=0.
$$
\end{prop}
and
\begin{prop}\label{lemma_5.6}
If the map $f_0$ for $e_0$ is singular, then
$$
\lim_{u\to 0}I_{F_T}(\lambda,u)=0.
$$
\end{prop}
Clearly, Proposition \ref{lemma_5.5} and  \ref{lemma_5.6} imply
Theorem \ref{prop_5.3}.

\v\n Suppose that $T$
is a graph. Let $C$ be a curve in $F_T$.

First assume that all components in $C$ are
smooth. Then copied from \cite{GP}, we have
$$
0\to \mc O_C\to\bigoplus_{\mbox{vertices}}\mc O_{C_v} \oplus
\bigoplus_{\mbox{edges}}\mc O_{Ce}\to \bigoplus_{\mbox{flags}} \mc
O_F\to 0,
$$
then (write $E=f^\ast T M^s_r$)
\begin{eqnarray*}
&&0 \to H^0(C,E) \to \bigoplus_{\mbox{vertices}} H^0(C_v, E)
\oplus \bigoplus_{\mbox{edges}}H^0(C_e,E)\\
&&\to \bigoplus_{\mbox{flags}} E_{p(F)}
\to H^1(C,E)\\
&& \to\bigoplus_{\mbox{edges}} H^1(C_v, E) \oplus
\bigoplus_{\mbox{edges}}H^1(C_e,E) \to 0.
\end{eqnarray*}
Please refer to \cite{GP} for flags. Here, by $p(F)$ we mean the
fixed point assigned to the flag. Hence the contribution of
$H^1/H^0$ is
\begin{equation}\label{eqn_new_6}
\frac{H^1(C,E)}{H^0(C,E)}
=\frac{\bigoplus_{\mbox{vertices}}E_{p(v)}^{val(v)-1}
\oplus\bigoplus_{\mbox{vertices}}H^1(C_e,E)}
{\bigoplus_{\mbox{edges}}H^0(C_e,E)}
\end{equation}
We translate each term to the equivariant Euler class, i.e, a
polynomial in $\lambda$ and $u$.

If $C$ is not smooth, each space in the long complex
should be replaced by the
invariant subspace with respect to
the proper finite group actions. Hence, each term in
the right hand side of \eqref{eqn_new_6} should be replaced
accordingly.

Recall that
\begin{equation}\label{eqn_new_7}
I^\Omega_{F_T}(\lambda,u)=
\int_{F_T}\frac{\Omega_{T^2}}{e_{T^2}(N^{vir}_{F_T})}.
\end{equation}
It is known that the equivariant form of $H^1/H^0$
gives
\begin{equation}\label{eqn_new_8}
\frac{1}{e_{T^2}(N^{vir}_{F_T})}.
\end{equation}
It is easy to show that
for each $e$, whose  map $f_e$  does not $\Gamma^s_r$,
$$
e_{T_2}\left(\frac{H^1(C_e,f_e^\ast E)}{H^0(C_e,f_e^\ast E)}
\right)
$$
contains no either $u$ or $u\inv$. We now focus other terms
in $H^1/H^0$.

\v\n {\em Claim 1:  if $f_0$ is smooth,
the equivariant  Euler form
 of the above term $H^1/H^0$
contains a positive power of $u$.}
 \v\n We count the possible
contributions for $u$. Case 1, there is $v=\mfp$ or $\mfq$ with
$val(v)>1$, we then have $u^{val(v)-1}$; Case 2, there is a
component $C_e$ that is a multiple of $\Gamma^s_r$, we then
actually have $u$ in the denorminator for $H^0(C_e,f_0^\ast\mc O)$ and
$ru$ in the numerator for $H^1(C_e,f_0^\ast\mc O(-2))$.  So they cancel
out.  Hence we have the claim. $\Box$

\v\n {\em Claim 2: if $f_0$ is an orbifold map,
  the equivariant Euler of
$H^1/H^0$ contains a factor of positive power of $u$.} \v\n

\v\n
Suppose $f_0$ is given by
$$
f_0: [S^2]\to \Gamma^s_r\subset W^s_r.
$$
Such a map can be realized by a map
$$
\tilde f_0: S^2 \to \tilde\Gamma^s_r\subset \tilde W^s_r
$$
with a quotient by $\mathbb Z_r$.
Here $[S^2]=S^2/\mathbb Z_r$.

Unlike the smooth case, $H^0([S^2],f_0^\ast E)$
contains no $u$. By the computations given below, in Corollary
\ref{cor_5.1} we
conclude that $H^1([S^2],f_0^\ast E)$ contains a factor $u$.
$\Box$

\v\n
We compute $H^1([S^2],f_0^\ast E)$
and its weight.
\def \inv{^{-1}}

Suppose that $\tilde f_0$ is a $d$-cover. Then on the sphere $\tilde
S^2$, the torus action weight at $0$ is $\lambda/d$ and at
$\infty$ is $-\lambda/d$,
and suppose that  $\mathbb Z_r$ action is $\mu$ at 0
and $\mu\inv$ at $\infty$;  for the pull-back bundle $\mc O(-2d)$
the torus action weight at fiber over $0$ is $-\lambda+ru$ and at
fiber over $\infty$ is $\lambda+ru$, the $\mathbb Z_r$ action are
$\mu^{-d}$ and $\mu^{d}$ for  the fibers over at $0$ and $\infty$.
These data are ready for us to compute the action on
$H^1([S^2],f_0^\ast E)$.

By Serre-duality, we have
$$
H^1(S^2, \mc O(-2d))=(H^0(S^2, \mc O(2d-2)))^\ast.
$$
The induced torus action weights on $\mc O(2d-2)$ at the fibers
over $0$ and $\infty$ are
$$
\frac{d-1}{d}\lambda-ru, \;\;\; \frac{1-d}{d}\lambda-ru
$$
respectively. The induced $\mathbb Z_r$ action are $\mu^{d-1}$ and
$\mu^{-d+1}$ for fibers on $\mc O(2d-2)$ over $0$ and $\infty$.

\begin{lemma}\label{lemma_5.7}
The sections of $H^0(S^2, \mc O(2d-2))$ are given by
$$
\{x^ay^b| a+b=2d-2, a,b\geq 0\}.
$$
The torus action weight for section $x^ay^b$ is
$\frac{d-1-a}{d}\lambda -ru$. The action  of $\mu$ is
$\mu^{d-1-a}$.
\end{lemma}
Hence
\begin{lemma}\label{lemma_5.8}
The section $x^ay^b$ that is $\mathbb Z_r$-invariant if and only
if $r|d-1-a$, and the torus action weight is
$\frac{d-1-a}{d}\lambda -ru$.
\end{lemma}
\begin{corollary}\label{cor_5.1}
$H^1([S^2], f_0^\ast\mc O(-2))$ contains a factor  $ru$.
\end{corollary}
{\bf Proof. } By the above lemma, we know that $x^{d-1}y^{d-1}$ is
$Z_r$-invariant and it is action weight is $-ru$. By taking the
dual, the corresponding factor is $ru$. q.e.d.

\v
\v\n
{\bf Proof of Proposition \ref{lemma_5.5}}: Claim 1 implies the proposition.
\v\n
{\bf Proof of Proposition \ref{lemma_5.6}}: Claim 2 implies the proposition.

\subsection{Vanishing results on $I^\Omega_{F_T}(\lambda,u)$, (II)}
\v\n Now suppose that $(\Gamma,\mc T,\mba)$ is an admissible data
given in Proposition \ref{prop_5.2}.
\begin{theorem}\label{thm_5.1}
If $(\Gamma,\mc T,\mba)$ is of case (1) and (2) in Proposition
\ref{prop_5.2},
$$
\langle\mba|\mbb\rangle_{\Gamma,\mc T}^{(M^s,Z)}=0.
$$
\end{theorem}
{\bf Proof. }Suppose that
$$
\langle\mba|\mbb\rangle_{\Gamma,\mc T}^{(M^s,Z)}=I(\Omega)
$$
for some $\Omega$ which has equivariant extension $\Omega_{T^2}$.
(cf. Lemma \ref{LEMMA_5.1}).

Let $F_T$ be a fix component of the torus action. It contributes 0
to the invariant unless the fixed curves $C\in F_T$ contains no
component covering $\Gamma^s$ (cf. Theorem \ref{prop_5.3}). It is easy to
conclude that $C$ {\em must} contains a ghost map
$$
f: (S^2, q_1,q_2,\cdots, q_l)\to M^s
$$
such that the image of $f$ is either $\mfp$ or $\mfq$ and sum of
the degree shifting numbers for all twisted sectors defined by
$q_i$ is $2+ l$. Suppose $\mfp=f(S^2)$. Again, we claim that
$I^\Omega_{F_T}(\lambda, u)$ contains  factor $u$. In fact,
$$
e_{T^2}(H^1(S^2,f^\ast \mc O_\mfp))= u.
$$
This is proved in \cite{CH}. Hence, it is easy to see the claim
follows. q.e.d.

\begin{defn}\label{defn_5.2}
If $$ \lim_{u\to 0}I^\Omega_{F_T}(\lambda,u)=0,
$$
we say the component $F_T$ contributes trivial to the invariant
$I(\Omega)$.
\end{defn}

\begin{corollary}\label{cor_5.2}
Let $(\Gamma,\mc T,\mba)$ be admissible, and
$$
I(\Omega)=\langle\mba|\mbb\rangle_{\Gamma,\mc T}^{(M^s,Z)}.
$$
If $F_T$ contributes nontrivial to the invariant, then for any
curve $C\in F_T$, all its connected components in the root
component $M^s$  must cover $L^s_{p,y}$ or $L^s_{q,x}$.
\end{corollary}
{\bf Proof. } By the previous theorem, the admissible data must be
of the 3rd case in Proposition \ref{prop_5.2}. Hence, points in
$C\cap Z$ must be either $\mfx$ or $\mfy$. Therefore, the
invariant curves that $C$ lives on must be $L^s_{p,y}$ and
$L^s_{q,x}$. q.e.d.

\section{Proof of the Main theorem}\label{section_6}
Combining \S\ref{section_2.3} and Theorem \ref{theorem_3.1}, we reduce Theorem \ref{thm_0.1} to
\begin{theorem}\label{thm_6.2}
$$
\sum_{A\not\in\mathbb
Z\Gamma^s_r}\langle\alpha_1^s,\alpha_2^s,\alpha_3^s\rangle_A
=\sum_{A\not\in\mathbb
Z\Gamma^{sf}_r}\langle\Psi^\ast\alpha_1^s,\Psi^\ast\alpha_2^s,\Psi^\ast\alpha_3^s\rangle_A.
$$
\end{theorem}
The rest of the section is devoted to the proof of this theorem.

\begin{remark}
We should point out that "$=$" is rather strong from the point of view of Ruan's conjecture.
Usually, it is conjectured that
$$
\langle\alpha_1^s,\alpha_2^s,\alpha_3^s\rangle\cong \langle\Psi^\ast\alpha_1^s,\Psi^\ast\alpha_2^s,\Psi^\ast\alpha_3^s\rangle.
$$
By "$\cong$", we mean that both sides equal up to analytic continuations. This is necessary
when
classes $[\Gamma^s]$ and $[\Gamma^{sf}]$ involved. For example, in \cite{CLZZ} we proved
$$
\langle\alpha_1^s,\alpha_2^s,\alpha_3^s\rangle_R
\cong \langle\Psi^\ast\alpha_1^s,\Psi^\ast\alpha_2^s,\Psi^\ast\alpha_3^s\rangle_R
$$
for Theorem \ref{theorem_3.1}.

But for the invariants that correspond to $A\not= d[\Gamma^s]$, it turns out
that we do not need the analytic continuation argument, the reason is because of
Theorem \ref{prop_5.3}: as long as $[\Gamma^s]$ (so is $[\Gamma^{sf}]$)
appears, the invariant vanishes.
\end{remark}

\subsection{Reducing the comparison to local models}\label{sect_6.2}
We now apply the degeneration formula to reduce the comparing
three point functions only  on local models. We explain this:
consider a three point function
$$
\langle\alpha_1^s,\alpha_2^s,\alpha_3^s\rangle_{0,A_s}.
$$
First, we observe that $[A_s]=A_s$
 since $\pi_\ast$ has no kernel (cf
 \cite{LR}). Denote the topology data by
$$
\Gamma^s=(0, A_s)
$$
and forms by $\mba^s=(\alpha_1^s,\alpha^s_2,\alpha^s_3)$.
Correspondingly, on $Y^{sf}$  we introduce
\begin{eqnarray*}
&& \mba^{sf}=(\Psi^\ast)\inv\mba^s,\\
&& A_{sf}=\Psi_\ast A_s,\\
&&\Gamma^{sf}=(0, A_{sf}).
\end{eqnarray*}
We write $\Gamma^{sf}=\Psi(\Gamma^s)$.

Consider the degenerations
\begin{eqnarray*}
&&\pi_s:Y^s\xrightarrow{degenerate} M^s_r\cup_Z Y^-;\\
&&\pi_{sf}:Y^{sf}\xrightarrow{degenerate}  M^{sf}_r\cup_Z Y^-.
\end{eqnarray*}
Let $\eta^s=(\Gamma^{+,s},\Gamma^{-},I_\rho)$ be a possible
splitting of $\Gamma^s$. Correspondingly,
$$
\Psi(\eta^s):=(\Psi_r(\Gamma^{+,s}),\Gamma^{-},I_\rho)
$$
gives a splitting of $\Psi(\Gamma^s)$. On the other hand, suppose
that
$$
\mba^s=\pi_s^\ast(\mba^{+,s},\mba^{-}).
$$
Then by the diagram \eqref{eqn_4.3},
$$
\mba^{sf}:=\Psi^\ast(\mba^s)=
\pi_{sf}^\ast(\Psi^\ast_r(\mba^{+,s}), \mba^{-}).
$$

\begin{prop}\label{prop_6.1}
 Suppose that  $\mba^s$ and $\mbb$ are given on $M^s_r$, then
\begin{equation}\label{eqn_6.1}
\langle \mba^{+,s}|\mbb,\mc T\rangle^{\ast(M^s,Z)}_{\Gamma^{+,s}}
=\langle \Psi^\ast_r\mba^{+,s}|\mbb,\mc
T\rangle^{\ast(M^s,Z)}_{\Psi_r(\Gamma^{+,s})}.
\end{equation}
Unlike $\langle\rangle^{\bullet(M^s,Z)}$, here
$\langle\rangle^{\ast(M^s,Z)}$ only sums over all admissible data.
\end{prop}

\begin{prop}
Proposition \ref{prop_6.1} $\Rightarrow$ Theorem \ref{thm_6.2}.
\end{prop}
{\bf Proof. }Applying the degeneration formula to
$\langle\mba^s\rangle_{\Gamma^s}$, we have
$$
\langle \mba^s\rangle_{\Gamma^s}^{Y^s}
=\sum_{I}\sum_{\eta=(\Gamma^{+,s},\Gamma^{-},I_\rho)}
C_\eta\langle \mba^{+,s}|\mbb^I,\mc
T\rangle^{\bullet(M^s,Z)}_{\Gamma^{+,s}} \langle
\mba^{-}|\mbb_I,\mc T\rangle^{\bullet(Y^{-},Z)}_{\Gamma^-}.
$$
Similarly,
$$
\langle \mba^{sf}\rangle_{\Gamma^{sf}}^{Y^{sf}}
=\sum_{I}\sum_{\eta^{sf}=(\Gamma^{+,sf},\Gamma^{-},I_\rho)}
C_\eta\langle \mba^{+,sf}|\mbb^I,\mc
T\rangle^{\bullet(M^s,Z)}_{\Gamma^{+,sf}} \langle
\mba^{-}|\mbb_I,\mc T\rangle^{\bullet(Y^{-},Z)}_{\Gamma^-}.
$$
Here $\ast^{sf}$ is always the correspondence of $\ast^s$ via
$\Psi$ or $\Psi_r$.

Since only admissible data contributes on the right hand sides of
two equations, \eqref{eqn_6.1} implies
$$
\langle \mba^s\rangle_{\Gamma^s}^{Y^s}=\langle
\mba^{sf}\rangle_{\Gamma^{sf}}^{Y^{sf}}.
$$
which is exactly what Theorem \ref{thm_6.2} asserts.

\subsection{Proof of Proposition \ref{prop_6.1}}\label{sect_6.3}
We now proceed to proof Proposition \ref{prop_6.1}.  Since the
moduli spaces in local models admit torus actions. By
localizations, we know the contributions only come from those fix
loci. So it is sufficient to compare fix loci and the invariants
they contribute.

For  $(M^s,Z)$, let $T$ be a graph, and $F_T^s$ be the component
of fix loci. Then by Corollary \ref{cor_5.2}, $F_T^s$ makes a
nontrivial contribution only when
 each curves in $F_T$ consists of only components on
$L^s_{p,y}$ and $L^s_{q,x}$.

Similarly, for $(M^{sf},Z)$, $F_T^{sf}$ makes a  nontrivial
contribution only when
 each curves in $F_T^{sf}$ consists of only components on
$L^{sf}_{p,y}$ and $L^{sf}_{q,x}$.

Since, the flop identifies
$$
L^s_{p,y}\leftrightarrow L^{sf}_{p,y}, L^s_{q,x}\leftrightarrow
L^{sf}_{q,x}
$$
and their normal bundles, therefore the flop identifies $F_T^s$
and $F_T^{sf}$ and their virtual normal bundles in their moduli
spaces. Hence,
$$
I_{F_T^s}(\lambda,u)=I_{F_T^{sf}}(\lambda,u).
$$
Proposition \ref{prop_6.1} then follows.

\end{document}